\documentclass[11pt]{article}
\usepackage{mathrsfs}
\usepackage{amsfonts}
\usepackage{dsfont}
\usepackage[notref,notcite]{showkeys}
\usepackage{CJK}
\usepackage{amsfonts,amssymb,amsmath,mathrsfs}

\usepackage{color,latexsym,amsfonts}
 \newtheorem{theorem}{Theorem}[section]
 \newtheorem{lemma}[theorem]{Lemma}
 \newtheorem{corollary}[theorem]{Corollary}
 \newtheorem{proposition}[theorem]{Proposition}
 
 \newtheorem{Definition}[theorem]{Definition}
 \newtheorem{remark}[theorem]{Remark}
 \newtheorem{condition}[theorem]{Condition}

 \def\blemma{\begin{lemma}\sl{}\def\elemma{\end{lemma}}}
 \def\btheorem{\begin{theorem}\sl{}\def\etheorem{\end{theorem}}}
 \def\bcorollary{\begin{corollary}\sl{}\def\ecorollary{\end{corollary}}}

 \def\bremark{\begin{remark}\sl{}\def\eremark{\end{remark}}}

 \def\beqlb{\begin{eqnarray}}\def\eeqlb{\end{eqnarray}}
 \def\beqnn{\begin{eqnarray*}}\def\eeqnn{\end{eqnarray*}}

 \def\<{\langle}\def\>{\rangle}

 \def\eqref#1{{\rm(\ref{#1})}}

\newcommand{\ra}{\rightarrow}

\def\d{\textup{d}}

\def\e{\textup{e}}

\def\fin{\hfill$\square$}

\def\newdot{{\kern.8pt\cdot\kern.8pt}}

\def\R{\mathbb{R}}
\def\E{\mathbb{E}}
\def\P{\mathbb{P}}

\def\F{\mathcal{F}}

\def\<{\langle}
\def\>{\rangle}
\def\Proof.{\noindent{\bf Proof.}}

\def\al{\alpha}
\def\be{\beta}
\def\ga{\gamma}
\def\de{\delta}
\def\ep{\epsilon}

\def\et{\eta}
\def\th{\theta}

\def\la{\lambda}
\def\rh{\rho}

\def\si{\sigma}
\def\ta{\tau}
\def\ph{\phi}

\def\om{\omega}

\def\ff{\frac}

\def\1{\mathds{1}}

\allowdisplaybreaks

\numberwithin{equation}{section}

\begin{document}

\

\noindent{}

\bigskip\bigskip

\centerline{\Large\bf Moment estimates and applications for SDEs driven by}

\smallskip

\centerline{\Large\bf fractional Brownian motion with irregular drifts}

\smallskip
\
\bigskip\bigskip

\centerline{Xi-Liang Fan$^\dag$, Shao-Qin Zhang$^\ddag$\footnote{E-mail address: fanxiliang0515@163.com (X. Fan), zhangsq@cufe.edu.cn(S. Zhang).}}

\bigskip

\centerline{\footnotesize{$^\dag$School of Mathematics and Statistics, Anhui Normal University, Wuhu 241002, China}}
\centerline{\footnotesize{$^\ddag$School of Statistics and Mathematics, Central University of Finance and Economics, Beijing 100081, China}}

\smallskip

\bigskip\bigskip

{\narrower{\narrower

\noindent{\bf Abstract.} In this paper, high-order moment, even exponential moment, estimates are established for the H\"older norm of solutions to stochastic differential equations  driven by fractional Brownian motion whose drifts are measurable and have linear growth. As applications, we first study the weak uniqueness of solutions to fractional stochastic differential equations. Moreover,  combining our estimates and the Fourier transform, we establish the existence of density of solutions to equations with irregular drifts. }


\bigskip



\section{Introduction}

The fractional Brownian motion appears naturally in the modeling of many complex phenomena in applications
when systems are subject to the rough external forcing.
The properties that the fractional Brownian motion with Hurst parameter $H\neq1/2$ is neither a Markov process nor a (weak) semimartingale
complicate the study of stochastic differential equations (SDEs in short) driven by fractional Brownian motion,
and then new techniques  beyond the classical Markovian framework are needed to construct such equations.
Based on the techniques of fractional calculus developed by Z\"{a}hle in \cite{Zahle98a},
Nualart and R\u{a}\c{s}canu \cite{Nualart&Rascanu02a} proved the existence and uniqueness result with $H>1/2$.
Coutin and Qian \cite{Coutin&Qian02a} also derived the existence and uniqueness result for $H\in(1/4,1/2)$ through
the rough-type arguments introduced in \cite{Lyons98a}.
For other results on the existence and uniqueness,
the reader may consult e.g. \cite{Boufoussi&Hajji11,Ferrante&Rovira06,Hu&Nualart&Songa,Liu&Yan18a,Neuenkirch&Nourdin&Tindel08,Nualart&Ouknine02b} and the references therein.

SDEs with irregular coefficients are intensively studied recently, see e.g. \cite{BM,KR,ZhangX16} and references therein. In this paper, we shall investigate SDEs driven by fractional Brownian motion on $\R^d$ of the form:
\begin{align}\label{equ_3_1}
\d X(t)=b(X(t))\d t+\sigma(X(t))\d B_t^H,\ \ X(0)=x, \ \ t\in[0,T],
\end{align}
with $H\in(1/2,1)$ under the hypothesis on
$$b:\R^d\rightarrow\R^d,\qquad \si:\R^d\ra \R^d\otimes\R^d$$
as follows:

\noindent{\bf {Hypothesis}} (A):
\begin{itemize}
\item[(i)] $b$ is measurable and has  linear growth;
\par
\item[(ii)] $\sigma$ is H\"{o}lder continuous of order $\gamma\in(1/H-1,1]$.
\end{itemize}

Moment estimates of solutions  play basic roles in studies of stochastic equations, and the exponential integrability of solutons is also of great use, see e.g. \cite{W17}.  In \cite{Duncan&Nualart09a}, the authors established some moment estimates to the solutions of \eqref{equ_3_1} under the assumptions that $b,\si$ are bounded and $\si$ is Lipschitz. \cite{Hu&Nualart07a} investigated moment estimates and  the exponential integrability of solutions to equations without drift term.  However, only  the uniform norm of solutions was considered in \cite{Hu&Nualart07a,Duncan&Nualart09a}.   Moment estimates and the exponential integrability of solutions to \eqref{equ_3_1}, especially for the associated H\"older norm,  under the mild condition (A) is also not known.   Our first result is to establish  moments estimates and the exponential integrability of solutions to fractional SDE \eqref{equ_3_1}, see Theorem \ref{Theorem 2.1}. As applications, we investigate the Girsanov transform associated with \eqref{equ_3_1} and obtain the uniqueness in law of the solutions to  \eqref{equ_3_1}.

The second main result of the present research is to establish the existence of the density for the SDE driven by  fractional Brownian motion satisfying the hypothesis (A). For studying the existence and the smoothness of densities of solutions to SDEs, Malliavin calculus is a powerful tool. In \cite{Cass&Friz10a,Cass&Friz&Victoir09a,Hu&Nualart&Songa,Nourdin&Simon06,Nualart&Saussereau09},
by rough path techniques or Malliavin calculus the authors handled the existence of densities of the solutions.  However, one usually requires regularity on the coefficients of the associated SDEs when using the Malliavin calculus.
There are some attempts to deal with irregularities in the coefficients of SDEs  without  using Malliavin calculus.
For instance, in \cite{Kohatsu-Higa&Makhlouf13a,Kohatsu-Higa&Tanaka12a} the authors obtained regularity properties
and upper and lower bounds for the densities of functionals of SDEs whose drifts are bounded and measurable,
through the Girsanov transform combined with an It\^{o}-Taylor expansion of the change of measure. On the other hand,  Fournier and Printems  introduced  a  method in \cite{Fournier&Printems10a} to investigate the existence of  densities  of some one-dimensional processes,
which do not have Malliavin derivatives including e.g. SDEs with H\"{o}lder coefficients and SDEs with random coefficients driven by usual Brownian motion. They used the characteristic function of random variable (i.e. the Fourier transform) and the one-step Euler approximation of the underlying processes.  This is developed and improved by Romito \cite{Ro18}, which can be applied to multidimensional stochastic equations with rough coefficients. However, Romito commented that this method might rule out processes driven by fractional Brownian motion, see the comments at the beginning of \cite[Section 7]{Ro18}. In this paper, we generalize the method used in \cite{Fournier&Printems10a,Ro18} to equations with fractional Brownian motion as the noise by characterizing the conditionally characteristic function of the one-step Euler scheme
and related estimates.  It is nontrivial.  For more details, see Lemma \ref{Lemma 3.l} and the proof of Theorem  \ref{Theorem 3.1} below.

This paper is structured as follows.
Section 2 is devoted to recall some useful facts on fractional calculus and fractional Brownian motion $B^H$. In Section 3,  we shall  present our results on the moments estimates and exponential integrability of solutions to \eqref{equ_3_1}. Some direct consequences, such as weak uniqueness of solutions, will be given there. In Section 4 we show that the law of $X(t)$ admits a density in some Besov space   for every $t>0$ under the hypothesis (A). 

\section{Preliminaries}

\subsection{Fractional calculus}

In this part, we shall give a brief account on fractional operators, which can be found in \cite{Samko&Kilbas&Marichev}.

Let $a,b\in\R$ with $a<b$.
For $\alpha>0$ and $f\in L^1(a,b)$,
the left-sided (resp. right-sided) fractional Riemann-Liouville integral of $f$ of order $\alpha$ on $[a,b]$ is defined as
\beqnn
I_{a+}^\alpha f(x)=\frac{1}{\Gamma(\alpha)}\int_a^x\frac{f(y)}{(x-y)^{1-\alpha}}\d y,\ \
\left(\mbox{resp.}~I_{b-}^\alpha f(x)=\frac{(-1)^{-\alpha}}{\Gamma(\alpha)}\int_x^b\frac{f(y)}{(y-x)^{1-\alpha}}\d y\right),
\eeqnn
where $x\in(a,b)$ a.e., $(-1)^{-\alpha}=\e^{-i\alpha\pi},\Gamma$ stands for the Euler function.
In particular, when $\alpha=n\in\mathbb{N}$, they reduced to the usual $n$-order iterated integrals.

Let $p\geq1$.
If $f\in I_{a+}^\alpha(L^p)$ (resp. $I_{b-}^\alpha(L^p)$) and $0<\alpha<1$,
then the Weyl derivatives read as follow
\beqlb\label{2.0'}
D_{a+}^\alpha f(x)=\frac{1}{\Gamma(1-\alpha)}\left(\frac{f(x)}{(x-a)^\alpha}+\alpha\int_a^x\frac{f(x)-f(y)}{(x-y)^{\alpha+1}}\d y\right)
\eeqlb
and
\beqlb\label{2.0''}
D_{b-}^\alpha f(x)=\frac{(-1)^\alpha}{\Gamma(1-\alpha)}\left(\frac{f(x)}{(b-x)^\alpha}+\alpha\int_x^b\frac{f(x)-f(y)}{(y-x)^{\alpha+1}}\d y\right),
\eeqlb
where the convergence of the integrals at the singularity $y=x$ holds pointwise for almost all $x$ if $p=1$ and in $L^p$-sense if $1<p<\infty$.

For any $\alpha\in(0,1)$, let $C^\alpha(a,b)$ be the space of $\alpha$-H\"{o}lder continuous functions $f$ on the interval $[a,b]$ and set
\beqnn
\|f\|_{a,b,\alpha}=\sup\limits_{a\leq s<t\leq b}\frac{|f(t)-f(s)|}{|t-s|^\alpha}.
\eeqnn
Besides, for any given continuous function $f:[a,b]\rightarrow\R$, put
\beqnn
\|f\|_{a,b,\infty}=\sup\limits_{a\leq s\leq b}|f(s)|.
\eeqnn

Suppose that $f\in C^\lambda(a,b)$ and $g\in C^\mu(a,b)$ with $\lambda+\mu>1$.
By \cite{Young36a}, the Riemann-Stieltjes integral $\int_a^bf\d g$ exists.
In \cite{Zahle98a}, Z\"{a}hle provides an explicit expression for the integral $\int_a^bf\d g$ in terms of fractional derivative. Let $\lambda>\alpha$ and $\mu>1-\alpha$.
Then the Riemann-Stieltjes integral can be expressed as
\beqlb\label{2.1}
\int_a^bf\d g=(-1)^\alpha\int_a^bD_{a+}^\alpha f(t)D_{b-}^{1-\alpha}g_{b-}(t)\d t,
\eeqlb
where $g_{b-}(t)=g(t)-g(b)$. \\
The relation \eqref{2.1} can be regarded as fractional integration by parts formula.

\subsection{Fractional Brownian motion}

For later use, we will recall some basic facts about fractional Brownian motion.
For a deeper discussion, we refer the reader to \cite{Biagini&Hu08a,Decreusefond&Ustunel98a,Mishura08a,Nualart06a}.

Let $B^H=\{B_t^H, t\in[0,T]\}$ be a fractional Brownian motion with Hurst parameter $H\in(1/2,1)$ defined on the probability
space $(\Omega,\F,\mathbb{P})$,
that is, $B^H$ is a centered Gaussian process with the covariance function
\beqnn
\E\left(B_t^HB_s^H\right)= R_H(t,s)=\frac{1}{2}\left(t^{2H}+s^{2H}-|t-s|^{2H}\right).
\eeqnn
Furthermore, one can show that $\E|B_t^H-B_s^H|^p=C(p)|t-s|^{pH},\ \forall p\geq 1$.
Consequently, by the Kolmogorov continuity criterion
$B^H$ have $(H-\epsilon)$-order H\"{o}lder continuous paths for all $\epsilon>0$.
For each $t\in[0,T]$, let $\mathcal {F}_t$ be the $\sigma$-algebra generated by the random variables $\{B_s^H:s\in[0,t]\}$ and the sets of
probability zero.

Denote $\mathscr{E}$ by the set of step functions on $[0,T]$.
Let $\mathcal {H}$ be the Hilbert space defined as the closure of
$\mathscr{E}$ with respect to the scalar product
\beqnn
\langle I_{[0,t]},I_{[0,s]}\rangle_\mathcal {H}=R_H(t,s).
\eeqnn
By B.L.T. theorem,
the mapping $I_{[0,t]}\mapsto B_t^H$ can be extended to an isometry between $\mathcal {H}$ and the Gaussian space $\mathcal {H}_1$ associated with $B^H$.
Denote this isometry by $\phi\mapsto B^H(\phi)$.

On the other hand, by \cite{Decreusefond&Ustunel98a} we know that the covariance kernel $R_H(t,s)$ can be written as
\beqlb\label{add2.1}
 R_H(t,s)=\int_0^{t\wedge s}K_H(t,r)K_H(s,r)\d r,
\eeqlb
where $K_H$ is a square integrable kernel given by
\begin{align}\label{equ-KH}
K_H(t,s)&=\ff {s^{1/2-H}} {\Gamma(H-1/2)}\int_s^t r^{H-1/2}(r-s)^{H-3/2}\d r\1_{[0,t]}(s)\nonumber\\
&=\Gamma\left(H+\frac{1}{2}\right)^{-1}(t-s)^{H-\frac{1}{2}}F\left(H-\frac{1}{2},\frac{1}{2}-H,H+\frac{1}{2},1-\frac{t}{s}\right),
\end{align}
in which $F(\cdot,\cdot,\cdot,\cdot)$ is the Gauss hypergeometric function (for details see \cite{Decreusefond&Ustunel98a} or \cite{Nikiforov&Uvarov88}).

Now, define the linear operator $K_H^*:\mathscr{E}\rightarrow L^2([0,T])$ by
\beqnn
(K_H^*\phi)(s)=K_H(T,s)\phi(s)+\int_s^T(\phi(r)-\phi(s))\frac{\partial K_H}{\partial r}(r,s)\d r.
\eeqnn
By integration by parts, it is easy to see that this can be rewritten as
\beqnn
(K_H^*\phi)(s)=\int_s^T\phi(r)\frac{\partial K_H}{\partial r}(r,s)\d r.
\eeqnn
Due to \cite{Alos&Mazet&Nualart01a}, for all $\phi,\psi\in\mathscr{E}$,
there holds $\langle K_H^*\phi,K_H^*\psi\rangle_{L^2([0,T])}=\langle\phi,\psi\rangle_\mathcal {H}$
and then $K_H^*$ can be extended to an isometry between $\mathcal{H}$ and $L^2([0,T])$.
Hence, according to \cite{Alos&Mazet&Nualart01a} again,
the process $\{W_t=B^H((K_H^*)^{-1}{\rm I}_{[0,t]}),t\in[0,T]\}$ is a Wiener process,
and $B^H$ has the following integral representation
\beqlb\label{add2.2}
 B^H_t=\int_0^tK_H(t,s)\d W_s.
\eeqlb

Finally, consider the operator $K_H: L^2([0,T])\rightarrow I_{0+}^{H+1/2 }(L^2([0,T]))$ associated with the integrable kernel $K_H(\cdot,\cdot)$
\beqnn
(K_Hf)(t)=\int_0^tK_H(t,s)f(s)\d s
\eeqnn
It can be proved (see \cite{Decreusefond&Ustunel98a}) that $K_H$ is an isomorphism and moreover,  for each $f\in L^2([0,T])$,
\beqnn
(K_H f)(s)=I_{0+}^{1}s^{H-1/2}I_{0+}^{H-1/2}s^{1/2-H}f.
\eeqnn
Consequently, for each $h\in I_{0+}^{H+1/2}(L^2([0,T]))$, the inverse operator $K_H^{-1}$ is of the form
\beqlb\label{2.2}
(K_H^{-1}h)(s)=s^{H-1/2}D_{0+}^{H-1/2}s^{1/2-H}h'.
\eeqlb

\section{Moment estimates of solutions to fractional SDE}

Fix any $T>0$.  In this section, we shall give some higher-order moment estimates of  solutions to \eqref{equ_3_1} on $[0,T]$ and present some interesting results induced by these estimates.  By \cite[Theorem 4]{Duncan&Nualart09a} and (A), it follows that there exists a solution $X$ to \eqref{equ_3_1}
which has $\beta$-H\"{o}lder continuous trajectories for any $\beta<H$. Generally, the condition (ii) ensures that the stochastic integral $\int_0^\cdot \si(X(t))\d B^H_t$ makes sense as a pathwise Riemann-Stieltjes integral. For a matrix $A\in\R^d\otimes\R^d$, we denote by $|A|$ the matrix norm of $A$, i.e.
$$|A|:=\sup_{|x|\leq 1, x\in\R^d}|Ax|.$$
In this and the following section, we shall denote by $\|f\|_\be$ and $\|f\|_\infty$ the $\be$-H\"older seminorm $\|f\|_{0,T,\be}$ and the supreme norm $\|f\|_{0,T,\infty}$ for simplicity.
Then our main result in this section reads as follows.
\btheorem\label{Theorem 2.1}
Assume (A).
Then   $\E\|X\|_{\beta}^p<\infty$  for all $p\geq 1$ and $0<\be<H$.
 Furthermore, if $\sigma$ is bounded, then for any $c>0$, $0<\de< 2H$ and $0<\be<H$,
it holds that $\E \exp[ c\|X\|_\be^{\de} ] <\infty.$
\etheorem
\bremark\label{Add Remark}
Let us compare with the relevant results proved in \cite[Theorem 4.1]{Hu&Nualart07a}.
They considered a SDE of the form:
\beqnn
X(t)=X(0)+\int_0^t\sigma(X(s))\d B_s^H, \ \ t\in[0,T],
\eeqnn
and proved the following integrability properties of $X$:
(i) If $\nabla\sigma$ is bounded and H\"{o}lder continuous of order larger than $\frac{1}{H}-1$,
then $\E\|X\|_\infty^p<\infty$ for all $p\geq2$;
(ii) If furthermore $\sigma$ is bounded, then $\E \exp[{ c\|X\|_\infty^{\de}}] <\infty$ holds for any $c>0$ and $0<\de< 2H$.

It is clear to see that our results apply to more general SDEs since we allow to have the drift term.
Moreover, we are able to present H\"{o}lder norm estimates under (A)
which requires only the H\"{o}lder continuity of $\sigma$.
\eremark

Before we prove this theorem, some interesting results induced by  these estimates will be presented.
As a direct consequence, we shall give some estimates on $R_t$ defined as follows
\beqnn
R_t:=\exp\left(\int_0^t\left\<K_H^{-1}\left(\int_0^\cdot h( X(r))\d r\right)(s),\d  W_s\right\>
-\frac{1}{2}\int_0^t\left|K_H^{-1}\left(\int_0^\cdot h( X(r))\d r\right)\right|^2(s)\d s\right),
\eeqnn
where $h$ is a $\R^d$-valued function on $\R^d$. These estimates may contribute to the study  of the exponential martingale of fractional Brownian motion.

\bcorollary\label{Cor2.1}
Assume (A).\\
(i) If $\si h$ has linear growth and
\begin{align}\label{h-Holder}
|h(x)-h(y)|\leq  K_1 |x-y|^\la(1+|x|^p+|y|^p),~~x,y\in\R^d,
\end{align}
with some  constants  $ K_1 >0$, $p>0$ and $\la\in (1-\ff 1 {2H},1]$. \\
Then $\{R_t\}_{t\in [0,T]}$ is a uniformly integrable martingale with
$$\sup_{t\in[0,T]}\E R_t\log R_t<\infty.$$
Consequently, $\displaystyle\left\{\tilde B^H_t\right\}_{t\in[0,T]}:=\left\{B^H_t-\int_0^t h(X(s))\d s\right\}_{t\in[0,T]}$  is a fractional Brownian motion under $R_T\P$. \\
(ii) If  $\si$ is bounded, and $h$ satisfies
\begin{align}\label{add-n-h}
|h(x)-h(y)|\leq  K_2 \left(|x-y|^{\la}\wedge 1\right)\left(1+|x|^q+|y|^q\right),~x,y\in\R^d
\end{align}
with some  constants  $ K_2 >0$,  $\la\in (1-\ff 1 {2H},1]$ and $q\in [0,H+\ff 1 {2H}-1)$,  then for any $ c >0$,
$$\E\exp\left( c \int_0^T\left|K^{-1}_H\left(\int_0^\cdot h(X(r))\d r\right)\right|^2(s)\d s\right)<\infty.$$
\ecorollary

As a direct consequence of Corollary \ref{Cor2.1}, we can obtain the uniqueness in law for weak solutions to equation \eqref{equ_3_1}.
Recall that by a weak solution to \eqref{equ_3_1} we mean a couple of adapted continuous processes $(B^H,X)$ on a filtered probability space
$(\Omega,\mathcal{F},\P,\{\mathcal{F}_t,t\in[0,T]\})$ such that
\begin{itemize}
\item[(i)] $B^H$ is an $\mathcal{F}_t$-fractional Brownian motion under $\P$,
\par
\item[(ii)] $(B^H,X)$ satisfy \eqref{equ_3_1}.
\end{itemize}

\bcorollary\label{Cor_3_2}
Assume (A).
Suppose that $\si$ is invertible, bounded, continuously differentiable and $\nabla \si$ is bounded and locally H\"older continuous of order $\ga'>\ff 1 H-1$ and
\begin{align}\label{h-Holder2}
|\si^{-1}b(x)-\si^{-1}b(y)|\leq L|x-y|^\la(1+|x|^p+|y|^p),~~x,y\in\R^d,
\end{align}
with some constants  $ L >0$, $p>0$ and $\la\in (1-\ff 1 {2H},1]$.
Then any two weak solutions to \eqref{equ_3_1} have the same distribution.
\ecorollary
Indeed, let $h=-\si^{-1}b$, then the conditions of Corollary \ref{Cor_3_2} imply that of Corollary \ref{Cor2.1}.
Thus, we derive that $\{R_t\}_{t\in [0,T]}$ is an exponential martingale.
So, the remaining part of the proof of Corollary \ref{Cor_3_2} can be treated thanks to the same kind of arguments of \cite[Theorem 8]{Duncan&Nualart09a},
and we omit the proof here.

In order to prove Theorem \ref{Theorem 2.1} and Corollary \ref{Cor2.1},
we begin with the following lemma which  play a key role in the proofs of their second assertions.
\blemma\label{BasIne_1}
Let $p>0$, $p'>1$ with $p-p'+1>0$. Then there exists a constant  $C(p,p')$ such that
\begin{align*}
\int_0^x\ff {u^{p-p'}} {1+k u^p}\d u\leq  C(p,p') \ff {x^{p-p'+1}} {1+x^{p-p'+1}k^{\ff {p-p'+1} p}},~~k>0, x>0.
\end{align*}
Consequently, it holds
$$\int_0^x\ff {k u^{p-p'}} {1+k u^p}\d u\leq  C(p,p') k^{\ff {p'-1} p},~~k>0, x>0.$$
\elemma

\emph{Proof.}
For $x>0$, we have
\begin{align*}
\int_0^x\ff {u^{p-p'}} {1+u^{p}}\d u &\leq \ff {x^{p-p'+1}} {p-p'+1}\1_{(0,1]}(x)+\left\{\ff 1 {p-p'+1}+\ff 1 {p'-1}\left(1-\ff 1 {x^{p'-1}}\right)\right\}\1_{(1,\infty)}(x)\\
&\leq \ff {p} {(p-p'+1)(p'-1)} \left(x^{p-p'+1}\wedge 1\right)\\
&\leq  C(p,p') \ff {x^{p-p'+1}} {1+x^{p-p'+1}},
\end{align*}
where in the last inequality we used the following inequality
\begin{align}\label{BasIne}
a\wedge 1\leq \ff {2a} {a+1},~~a\geq 0,
\end{align}
and $C(p,p')$ is a constant depending only on $p$ and $p'$.\\
Then, it follows that
\begin{align*}
\int_0^{x}\ff {u^{p-p'}} {1+ku^{p}}\d u&=k^{-\ff {p-p'+1} {p} }\int_0^{xk^{\ff 1 {p}} }\ff {v^{p-p'}} {1+v^{p}}\d v\\
&\leq  C(p,p') k^{-\ff {p-p'+1} {p}}\ff {x^{p-p'+1}k^{\ff {p-p'+1} {p}}} {1+x^{p-p'+1}k^{\ff {p-p'+1} {p}}}\\
&= C(p,p') \ff {x^{p-p'+1}} {1+x^{p-p'+1}k^{\ff {p-p'+1} {p}}}.
\end{align*}
Consequently, we have
\begin{align*}
\int_0^{x}\ff {ku^{p-p'}} {1+ku^{p}}\d u
&\leq  C(p,p') \ff {kx^{p-p'+1}} {1+x^{p-p'+1}k^{\ff {p-p'+1} {p}}}\\
&= C(p,p') \ff {k^{\ff {p-p'+1} {p}}x^{p-p'+1}} {1+x^{p-p'+1}k^{\ff {p-p'+1} {p}}}k^{\ff {p'-1} {p}}
\leq  C(p,p') k^{\ff {p'-1} {p}}.
\end{align*}
\fin

\bigskip

\emph{Proof of Theorem \ref{Theorem 2.1}.}
(1) We shall prove the first assertion.
 By \eqref{equ_3_1}, we deduce that, for any $0\leq s<t\leq T$,
\beqlb\label{Lemma 3.2-0}
X(t)-X(s)=\int_s^tb(X(r))\d r+\int_s^t\sigma(X(r))\d B_r^H.
\eeqlb
Since $\ga\in (\ff 1 H-1,1]$, $\ff {1} {1+\ga}<H$. Taking $H>\be>\ff 1 {1+\ga}$, then $\ga\be>1-\be$.
Using {\color{red} {the}} fractional by parts formula \eqref{2.1} with $\al\in (1-\be,\beta\gamma)$, we get
\beqlb\label{Lemma 3.2-1}
\int_s^t\sigma(X(r))\d B_r^H=(-1)^\alpha\int_s^tD_{s+}^\alpha\sigma(X(\cdot))(r)D_{t-}^{1-\alpha}B^H_{t-}(r)\d r.
\eeqlb
It follows from \eqref{2.0'} and \eqref{2.0''} that
\begin{align}\label{Lemma 3.2-2}
\left|D_{s+}^\alpha\sigma(X_\cdot)(r)\right|
= &\frac{1}{\Gamma(1-\alpha)}\left|\frac{\sigma(X(r))}{(r-s)^\alpha}+\alpha\int_s^r\frac{\sigma(X(r))-\sigma(X(u))}{(r-u)^{\alpha+1}}\d u\right|\nonumber\\
\leq &\frac{1}{\Gamma(1-\alpha)}\Big(|\sigma(X(r))|(r-s)^{-\alpha}\nonumber\\
&\qquad+\alpha\|\sigma\|_\gamma\int_s^r\|X\|_{u,r,\beta}^\gamma(r-u)^{\beta\gamma-\alpha-1}\d u\Big)\nonumber\\
\leq&C\left(|\sigma(X(r ))|(r-s)^{-\alpha}+\|\sigma\|_\gamma\|X\|_{s,r,\beta}^\gamma(r-s)^{\beta\gamma-\alpha}\right)
\end{align}
and
\beqlb\label{Lemma 3.2-3}
\left|D_{t-}^{1-\alpha}B^H_{t-}(r)\right|
&=&\frac{1}{\Gamma(\alpha)}\left|\frac{B^H_r-B^H_t}{(t-r)^{1-\alpha}}+(1-\alpha)\int_r^t\frac{B^H_r-B^H_u}{(u-r)^{2-\alpha}}\d u\right|\nonumber\\
&\leq&C\|B^H\|_\beta(t-r)^{\alpha+\beta-1},
\eeqlb
where and in what follows $C$ denotes a generic constant.\\
Then plugging \eqref{Lemma 3.2-2} and \eqref{Lemma 3.2-3} into \eqref{Lemma 3.2-1} yields
\beqlb\label{Lemma 3.2-4}
\left|\int_s^t\sigma(X(r))\d B_r^H\right|
&\leq&C\|B^H\|_\beta\Bigg(\int_s^t|\sigma(X(r))|(r-s)^{-\alpha}(t-r)^{\alpha+\beta-1}\d r\nonumber\\
&&\ \ \ \ \ \ \ \ \ \ \ +\|\sigma\|_\gamma\int_s^t\|X\|_{s,r,\beta}^\gamma(r-s)^{\beta\gamma-\alpha}(t-r)^{\alpha+\beta-1}\d r\Bigg).
\eeqlb
For the drift term, using the linear growth of $b$ we clearly get
\beqlb\label{Lemma 3.2-5}
\left|\int_s^tb(X(r))\d r\right|\leq C\int_s^t\left(1+|X(r)|\right)\d r.
\eeqlb
Then it follows from \eqref{Lemma 3.2-0}, \eqref{Lemma 3.2-4} and \eqref{Lemma 3.2-5} that
\begin{align}\label{fanadd-1}
\ff {|X(t)-X(s)|} {(t-s)^\be}\leq &C(t-s)^{1-\be}+\ff C {(t-s)^\be}\int_s^t|X(r)|\d r\nonumber\\
& +\ff {C\|B^H\|_\be} {(t-s)^\be}\int_s^t|\si(X(r))|(r-s)^{-\al}(t-r)^{\al+\be-1}\d r\nonumber\\
&+\ff {C\|\si\|_\ga\|B^H\|_\be} {(t-s)^\be}\int_s^t\|X\|_{s,r,\be}^\ga(r-s)^{\ga\be-\al}(t-r)^{\al+\be-1}\d r\nonumber\\
\leq &C(t-s)^{1-\be}(1+\|X\|_{s,t,\infty})\nonumber\\
& +\ff {C\|B^H\|_\be} {(t-s)^\be}\int_s^t|\si(X(s))|(r-s)^{-\al}(t-r)^{\al+\be-1}\d r\nonumber\\
& +\ff {C\|B^H\|_\be} {(t-s)^\be}\int_s^t|\si(X(r))-\si(X(s))|(r-s)^{-\al}(t-r)^{\al+\be-1}\d r\nonumber\\
&+\ff {C\|\si\|_\ga\|B^H\|_\be} {(t-s)^\be}\int_s^t\|X\|_{s,r,\be}^\ga(r-s)^{\ga\be-\al}(t-r)^{\al+\be-1}\d r\nonumber\\
\leq &C(t-s)^{1-\be}(1+\|X\|_{s,t,\infty})+C\|B^H\|_\be|\si(X(s))|\nonumber\\
&+\ff {C\|\si\|_\ga\|B^H\|_\be} {(t-s)^{\be}}\int_s^t\|X\|_{s,r,\be}^\ga(r-s)^{\ga\be-\al}(t-r)^{\al+\be-1}\d r\nonumber\\
\leq &C(t-s)^{1-\be}(1+\|X\|_{s,t,\infty})+C\|B^H\|_\be(1+|X(s)|^\ga)\nonumber\\
&+\ff {C\|\si\|_\ga\|B^H\|_\be} {(t-s)^{\be}}\int_s^t\|X\|_{s,r,\be}^\ga(r-s)^{\ga\be-\al}(t-r)^{\al+\be-1}\d r\nonumber\\
\leq & C\left((t-s)^{1-\be}+\|B^H\|_\be\right)(1+\|X\|_{s,t,\infty}) \nonumber\\
& +\ff {C\|\si\|_\ga\|B^H\|_\be} {(t-s)^{\be}}\int_s^t\|X\|_{s,r,\be}^\ga(r-s)^{\ga\be-\al}(t-r)^{\al+\be-1}\d r.
\end{align}

 Next, we shall focus on the upper bound of the right hand side in \eqref{fanadd-1}.

Since $\be>\ff 1 {1+\gamma}$ and $\ga\in(\ff 1 H-1,1]$, we get $\be>1-\ff 1 {2\ga}$.
Then,  there is $\th\in (0,1)$ such that $\ff 1 {2\ga}>\th>1-\be$, i.e.
$$\th+\be>1,\qquad \ff 1 {1-\ga\th}<2.$$
By  the Young inequality and  the H\"older inequality and $\ga\be-\al>0$, $\al+\be-1>0$, we have
\begin{align*}
&\ff {C\|\si\|_\ga\|B^H\|_\be} {(t-s)^{\be}}\int_s^t\|X\|_{s,r,\be}^\ga(r-s)^{\ga\be-\al}(t-r)^{\al+\be-1}\d r\\
\leq &\|X\|_{s,t,\be}^{\ga\th}\ff {C\|\si\|_\ga\|B^H\|_\be} {(t-s)^{\be}}\int_s^t\|X\|_{s,r,\be}^{\ga(1-\th)}(t-r)^{\al+\be-1}(r-s)^{\be\ga-\al}\d r\\
\leq &\ff 1 2\|X\|_{s,t,\be}+C\left(\ff {\|B^H\|_\be} {(t-s)^{\be}}\int_s^t\|X\|_{s,r,\be}^{\ga(1-\th)}(t-r)^{\al+\be-1}(r-s)^{\be\ga-\al}\d r\right)^{\ff 1 {1-\ga\th}}\\
\leq &\ff 1 2\|X\|_{s,t,\be}+C\|B^H\|_\be^{\ff 1 {1-\ga\th}}(t-s)^{\ff {\ga\be-1} {1-\ga\th}}\left(\int_s^t\|X\|_{s,r,\be}^{\ga(1-\th)}\d r\right)^{\ff 1 {1-\ga\th}}\\
\leq &\ff 1 2\|X\|_{s,t,\be}+C\|B^H\|_\be^{\ff 1 {1-\ga\th}}(t-s)^{\ff {\ga(\be+\th-1)} {1-\ga\th}}\left(\int_s^t\|X\|_{s,r,\be}\d r\right)^{\ff {\ga(1-\th)} {1-\ga\th}}\\
\leq &\ff 1 2\|X\|_{s,t,\be}+C\|B^H\|_\be^{\ff 1 {1-\ga\th}}(t-s)^{\ff {\ga(\be+\th-1)} {1-\ga\th}}\left(1+\int_s^t\|X\|_{s,r,\be}\d r\right).
\end{align*}
Hence, substituting this into \eqref{fanadd-1} leads to
\begin{align}\label{Xbeta}
 \ff {|X(t)-X(s)|} {(t-s)^\be}
\leq &C\left((t-s)^{1-\be}+\|B^H\|_\be\right)\left(1+\|X\|_{s,t,\infty}\right)\nonumber\\
& +C\|B^H\|_\be^{\ff 1 {1-\ga\th}}(t-s)^{\ff {\ga(\be+\th-1)} {1-\ga\th}}\nonumber\\
&+C\|B^H\|_\be^{\ff 1 {1-\ga\th}}(t-s)^{\ff {\ga(\be+\th-1)} {1-\ga\th}}\int_s^t\|X\|_{s,r,\be}\d r.
\end{align}
As for the estimate of the term $\|X\|_{s,t,\infty}$ in the right hand of \eqref{Xbeta},
by \eqref{Lemma 3.2-0}, \eqref{Lemma 3.2-4} and \eqref{Lemma 3.2-5} again we first obtain that
\begin{align*}
\left|X(t)\right|\leq &\left|X(s)\right|+C\int_s^t(1+|X(r)|)\d r+C\|B^H\|_\beta\int_s^t|\sigma(X(r))|(r-s)^{-\alpha}(t-r)^{\alpha+\beta-1}\d r\\
&+C\|B^H\|_\be\|\sigma\|_\gamma\int_s^t\|X\|_{s,r,\beta}^\gamma(r-s)^{\gamma\beta-\alpha}(t-r)^{\alpha+\beta-1}\d r\\
\leq & |X(s)|+C(t-s)(1+|X(s)|)+C\int_s^t\|X\|_{s,r,\beta}(r-s)^{\beta}\d r\\
&+C\|B^H\|_{\beta}|\sigma(X(s))|\int_s^t(r-s)^{-\alpha}(t-r)^{\alpha+\beta-1}\d r\\
&+C\|B^H\|_{\beta}\int_s^t|\si(X(r))-\sigma(X(s))|(r-s)^{-\alpha}(t-r)^{\alpha+\beta-1}\d r\\
&+C\|B^H\|_\beta\|\sigma\|_{\gamma}\int_s^t\|X\|_{s,r,\beta}^\gamma(r-s)^{\gamma\beta-\alpha}(t-r)^{\alpha+\beta-1}\d r\\
\leq & |X(s)|+C(t-s)(1+|X(s)|)+C(t-s)^\be\int_s^t\|X\|_{s,r,\beta}\d r\\
&+C\|B^H\|_{\beta}|\sigma(X(s))|(t-s)^\be\\
&+C\|B^H\|_\beta\|\sigma\|_{\gamma}\int_s^t\|X\|_{s,r,\beta}^\gamma(r-s)^{\gamma\beta-\alpha}(t-r)^{\alpha+\beta-1}\d r\\
\leq & |X(s)|+C(t-s)(1+|X(s)|)\\
&+C(t-s)^\be\int_s^t\|X\|_{s,r,\beta}\d r+C\|B^H\|_{\beta}(1+|X(s)|^\ga)(t-s)^\be\\
&+C\|B^H\|_\beta\|\sigma\|_{\gamma}\int_s^t( 1+\|X\|_{s,r,\beta})(r-s)^{\gamma\beta-\alpha}(t-r)^{\alpha+\beta-1}\d r\\
\leq &|X(s)|+C(t-s)(1+|X(s)|)\\
&+C\|B^H\|_\beta(1+|X(s)|^\gamma)(t-s)^\beta +C\|B^H\|_\beta(t-s)^{(1+\gamma)\beta}\\
&+C\left((t-s)^\be+(t-s)^{(1+\ga)\be-1}\|B^H\|_\be\right)\int_s^t\|X\|_{s,r,\be}\d r,
\end{align*}
where in the last inequality, we used the inequalities $\ga\be-\al>0$ and $\al+\be-1>0$. Since this inequality holds with $X(t)$ replaced by $X(t')$ for any $s\leq t'\leq t$,  we arrive at
\begin{align}\label{Xinfty}
\|X\|_{s,t,\infty}\leq &|X(s)|+C(t-s)(1+|X(s)|)\nonumber\\
&+C\|B^H\|_\beta(1+|X(s)|^\gamma)(t-s)^\beta+C\|B^H\|_\beta(t-s)^{(1+\gamma)\beta}\nonumber\\
&+C\left((t-s)^\be+(t-s)^{(1+\ga)\be-1}\|B^H\|_\be\right)\int_s^t\|X\|_{s,r,\be}\d r\nonumber\\
\leq &C(t-s)+C\left(1+(t-s)+(t-s)^\beta\|B^H\|_\be\right)|X(s)|\nonumber\\
&+C\left((t-s)^\be+(t-s)^{(1+\gamma)\beta}\right)\|B^H\|_\beta\nonumber\\
&+C\left((t-s)^\be+(t-s)^{(1+\ga)\be-1}\|B^H\|_\be\right)\int_s^t\|X\|_{s,r,\be}\d r.
\end{align}
Now, plugging \eqref{Xinfty} into \eqref{Xbeta} implies
\begin{align*}
\|X\|_{s,t,\be}\leq &C\left[(t-s)^{1-\be}+\|B^H\|_\be+\|B^H\|_\be^{\ff 1 {1-\ga\th}}(t-s)^{\ff {\ga(\be+\th-1)}{1-\ga\th}}\right]\\
&+C\left((t-s)^{1-\be}+\|B^H\|_\be\right)\Big\{(t-s)+\left(1+(t-s)+(t-s)^\be\|B^H\|_\be\right)|X(s)|\\
&\qquad+\left((t-s)^\be+(t-s)^{(1+\ga)\be}\right)\|B^H\|_\be\\
&\qquad+\left((t-s)^\be+(t-s)^{(1+\ga)\be-1}\|B^H\|_\be\right)\int_s^t\|X\|_{s,r,\be}\d r\Big\}\\
&+C\|B^H\|_\be^{\ff 1 {1-\th\ga}}(t-s)^{\ff {\ga(\be+\th-1)} {1-\ga\th}}\int_s^t\|X\|_{s,r,\be}\d r\\
\leq &C\left(1+(t-s)^{2-\be}\right)\left(1+\|B^H\|_\be\right)^2|X(s)|+C\left(1+(t-s)^{1+\ga\be}\right)\left(1+\|B^H\|_\be\right)^2\nonumber\\
&+C\Big[(t-s)+(t-s)^{\ff {\ga(\be+\th-1)} {1-\ga\th}}\|B^H\|_\be^{\ff 1 {1-\ga\th}}+(t-s)^{(1+\ga)\be-1}\|B^H\|_\be^2\nonumber\\
&\qquad+\|B^H\|_\be\left((t-s)^{\be}+(t-s)^{\ga\be}\right)\Big]\int_s^t\|X\|_{s,r,\be}\d r\nonumber\\
\leq &C\left(1+(t-s)^{1+\ga\be}\right)\left(1+\|B^H\|_\be\right)^2\Big(1+|X(s)|\Big)\nonumber\\
&+C\Big[(t-s)+(t-s)^{\ff {\ga(\be+\th-1)} {1-\ga\th}}\|B^H\|_\be^{\ff 1 {1-\ga\th}}+(t-s)^{(1+\ga)\be-1}\|B^H\|_\be^2\nonumber\\
&\qquad+\|B^H\|_\be\left((t-s)^{\be}+(t-s)^{\ga\be}\right)\Big]\int_s^t\|X\|_{s,r,\be}\d r,
\end{align*}
where in the last two inequalities we have used that $1+\ga\be>2-\be$.\\
Recalling that $\th+\be>1$ and $\ff 1 {1-\ga\th}<2$, by the Gronwall inequality we conclude that
\begin{align*}
\|X\|_{s,t,\be}\leq C\left(1+(t-s)^{1+\ga\be}\right)\left(1+\|B^H\|_\be\right)^2\Big(1+|X(s)|\Big)\exp\left(C(t-s)^{1+\de}(\|B^H\|_\be^2+1)\right)
\end{align*}
with $\delta=\min\left\{(1+\gamma)\beta-1,\ff {\ga(\be+\th-1)} {1-\ga\th}\right\}>0$.\\
Let $n\in\mathbb N$ and $k_0=t-s$.
Then, we get
\begin{align*}
&\|X\|_{(n-1)k_0,nk_0,\be}\\
\leq &C\left(1+k_0^{1+\ga\be}\right)\left(1+\|B^H\|_\be\right)^2\exp\left(Ck_0^{1+\de}(\|B^H\|_\be^2+1)\right)\Big(|X((n-1)k_0)|+1\Big)\\
\leq &Ck_0^\be\left(1+k_0^{1+\ga\be}\right)\left(1+\|B^H\|_\be\right)^2\exp\left(Ck_0^{1+\de}(\|B^H\|_\be^2+1)\right)\|X\|_{(n-2)k_0,(n-1)k_0,\be}\\
&+C\left(1+k_0^{1+\ga\be}\right)\left(1+\|B^H\|_\be\right)^2\exp\left(Ck_0^{1+\de}(\|B^H\|_\be^2+1)\right)\Big(|X((n-2)k_0)|+1\Big)\\
\leq &C\Big\{k_0^\be\left(1+k_0^{1+\ga\be}\right)^2\left(1+\|B^H\|_\be\right)^4\exp\left(2Ck_0^{1+\de}(\|B^H\|_\be^2+1)\right)\\
&+\left(1+k_0^{1+\ga\be}\right)\left(1+\|B^H\|_\be\right)^2\exp\left(Ck_0^{1+\de}(\|B^H\|_\be^2+1)\right)\Big\}\Big(|X((n-2)k_0)|+1\Big)\\
\leq &C (k_0^{\be}+1)\Big\{\left(1+k_0^{1+\ga\be}\right)\left(1+\|B^H\|_\be\right)^{2}\exp\left(Ck_0^{1+\de}(\|B^H\|_\be^2+1)\right)\Big\}^2\Big(|X((n-2)k_0)|+1\Big).
\end{align*}
We iterate to find the following inequality
\begin{align}\label{Fan-add9}
&\|X\|_{(n-1)k_0,nk_0,\be}\cr
\leq &C_n (k_0^{\be}+1)^{n-1}\Big\{\left(1+k_0^{1+\ga\be}\right)\left(1+\|B^H\|_\be\right)^{2}\exp\left(Ck_0^{1+\de}(\|B^H\|_\be^2+1)\right)\Big\}^n\Big(|X(0)|+1\Big)\cr
\end{align}
for some positive constant $C_n$. \\
By the Fernique theorem we know that there is $\et>0$ such that $\E\exp[\et\|B^H\|_\be^2]<\infty$. Thus for every $m\in\mathbb N$,
we let
$$k_0=\left(\ff {\et} {2Cnm}\right)^{\ff 1 {1+\de}},$$
which, together with \eqref{Fan-add9}, leads to $\E\|X\|_{(n-1)k_0,nk_0,\be}^m<\infty$ and then
$\E\|X\|_{0,nk_0,\be}^m<\infty.$
Since
$$\lim_{n\ra\infty} nk_0=\lim_{n\ra\infty}\left(\ff {\et } {2Cm}\right)^{\ff 1 {1+\de}}n^{\ff {\de} {1+\de}}=\infty,$$
we obtain that
$$\E\|X\|_\be^m<\infty,\qquad m\in\mathbb N$$
with $\be>\ff 1 {1+\ga}$  (then for all $\be<H$), which also implies that $\E\|X\|_\infty^m<\infty$ for all $m\in \mathbb N$.

\medskip

(2) Suppose that $\si$ is bounded in addition. Let $H>\be>\ff 1 {1+\ga}$. Then there is $\al$ such  that $1-\be<\al<\ga\be$. In view of \eqref{Lemma 3.2-0}-\eqref{Lemma 3.2-3} and the boundedness of $\si$, we derive
\begin{align}\label{3.1-2}
&|X(t)-X(s)|-\int_s^t|b(X(r))|\d r\leq \left|\int_s^t\si( X(r))\d  B_r^H\right|\nonumber\\
\leq &C\|B^H\|_\be\Big(\int_s^t|\si(  X(r))|(r-s)^{-\al}(t-r)^{\al+\be-1}\d r\nonumber\\
&\quad+\int_s^t\int_s^r\ff {|\si( X(r))-\si( X(u))|} {(r-u)^{\al+1}}\d u(t-r)^{\al+\be-1}\d r\Big)\nonumber\\
\leq &C \|  B^H\|_\be\Big(\|\si\|_\infty(t-s)^\be\nonumber\\
&\quad +(t-s)^{\al+\be-1}\int_s^t\int_s^r\ff {(\|\si\|_\ga\|  X\|_{u,r,\be}^\ga(r-u)^{\ga\be})\wedge (2\|\si\|_\infty)} {(r-u)^{\al+1}}\d u\d r\Big)\nonumber\\
\leq &C \|  B^H\|_\be\Big((t-s)^\be\nonumber\\
&\quad +(t-s)^{\al+\be-1}\int_s^t\int_s^r\ff {(\|  X\|_{u,r,\be}^\ga(r-u)^{\ga\be})\wedge 1} {(r-u)^{\al+1}}\d u\d r\Big)\nonumber\\
\leq & C\|  B^H\|_\be\left((t-s)^\be+(t-s)^{\al+\be-1}\int_s^t\int_s^r\ff {\| X\|_{u,r,\be}^\ga(r-u)^{\ga\be-\al-1}} {1+(r-u)^{\ga\be}\|  X\|_{u,r,\be}^\ga}\d u\d r\right)\nonumber\\
\leq &  C\|  B^H\|_\be\left((t-s)^\be+(t-s)^{\al+\be-1}\int_s^t\int_s^r\ff {\|  X\|_{s,r,\be}^\ga(r-u)^{\ga\be-\al-1}} {1+(r-u)^{\ga\be}\|  X\|_{s,r,\be}^\ga}\d u\d r\right).
\end{align}
Here, in the last  two inequalities, we resorted to \eqref{BasIne}
and the monotonicity of the function $\frac{x}{1+\tilde{c}x}$ with a constant $\tilde{c}$, respectively. \\
By Lemma \ref{BasIne_1}, we get
\begin{align*}
\int_s^r\ff {\|  X\|_{s,r,\be}^\ga(r-u)^{\ga\be-\al-1}} {1+(r-u)^{\ga\be}\| X\|_{s,r,\be}^\ga}\d u\leq &C\left( \| X\|_{s,r,\be}^\ga\right)^{\ff {\al} {\ga\be}}=C\|  X\|_{s,r,\be}^{\ff {\al} {\be}}.
\end{align*}
Substituting this into \eqref{3.1-2}, and taking into account of the condition $\mathrm{(i)}$ of (A), we have
\begin{align*}
|X(t)-X(s)|&\leq C \|  B^H\|_\be\left((t-s)^\be+(t-s)^{\al+\be}\| X\|_{s,t,\be}^{\ff {\al} {\be}}\right)+C\int_s^t(1+|X(r)|)\d r\nonumber\\
& \leq  C \|  B^H\|_\be\left((t-s)^\be+(t-s)^{\al+\be}\| X\|_{s,t,\be}^{\ff {\al} {\be}}\right)\nonumber\\
&\qquad +C(t-s)(1+|X(s)|)+C(t-s)^\be\int_s^t\|X\|_{s,r,\be}\d r.
\end{align*}
Then, the Gronwall inequality and the Young inequality imply
\begin{align}\label{3.1-add_1}
\|X\|_{s,t,\be}\leq Ce^{C(t-s)}\left[\| B^H\|_\be\left(1+(t-s)^{\be}\| X\|_{s,t,\be}\right)+(t-s)^{1-\be}(1+|X(s)|)\right].
\end{align}
Let $\Delta>0$ such that
\begin{align}\label{ad-De}
1-Ce^{CT}\Delta^\be \|B^H\|_\be =\ff 1 2.
\end{align}
Then, by taking $t=s+\Delta$  and \eqref{3.1-add_1}  we get
\begin{align}\label{ad-n-Xbe}
\|X\|_{s,s+\Delta,\be}&\leq \ff {Ce^{CT}\left[\|B^H\|_\be+\Delta^{1-\be}(1+|X(s)|)\right]} {1-Ce^{CT}\Delta^\be{ \|B^H\|_\be}}\nonumber\\
& =  2Ce^{CT}\left(\|B^H\|_\be+\Delta^{1-\be}(1+|X(s)|)\right).
\end{align}
Set
$$\Lambda_k=\|X\|_{k\Delta,(k+1)\Delta,\be}\Delta^\be,\quad \Gamma_k=|X(k\Delta)|,\quad \th=2Ce^{CT}\Delta,~k\in \mathbb{N}.$$
Then by \eqref{ad-n-Xbe} and \eqref{ad-De}, we obtain
\begin{align}\label{Fan-add1}
\Lambda_k&\leq 2Ce^{CT}\left(\|B^H\|_\be \Delta^\be+\Delta (1+|X(k\Delta)|)\right)\cr
&\leq  1+\th+\Gamma_k \th.
\end{align}
Noticing that
\begin{align*}
\Gamma_k & \leq \left|X(k\Delta)-X((k-1)\Delta)\right|+\left|X((k-1)\Delta)\right|\\
& \leq \Lambda_{k-1}+\Gamma_{k-1}\leq 1+\th+\Gamma_{k-1}(1+\th),
\end{align*}
we have
\begin{align*}
\Gamma_k\leq (1+\th)^k\Gamma_{0}+\sum_{j=1}^k(1+\th)^j=(1+\th)^k|X(0)|+\ff {(1+\th)^{k+1}-(1+\th)} {\th}.
\end{align*}
Hence,  combining this with \eqref{Fan-add1} yields
\begin{align*}
\Lambda_k\leq (1+\th)^{k+1}+\th(1+\th)^k|X(0)|\leq (1+\theta)^{k+1}(1+\theta |X(0)|).
\end{align*}
Then, we derive that
\begin{align}\label{Fan-add2}
\|X\|_\infty& \leq |X(0)|+\sum_{j=0}^{\left[\ff T {\Delta}\right]}\|X\|_{j\Delta,(j+1)\Delta,\be}\Delta^\be =|X(0)|+\sum_{j=0}^{\left[\ff T {\Delta}\right]}\Lambda_{j}\cr
&\leq |X(0)|+\ff {(1+\th)^{\left[\ff T {\Delta}\right]+2}-(1+\th)} {\th}\left(1+\th|X(0)|\right)\cr
& \leq(1+\th)^{ \ff T {\Delta} +2} |X(0)|+(1+\th)^{ \ff T {\Delta} +2}\th^{-1}\cr
& =(1+\th)^{ \ff T {\Delta} +2} \left(|X(0)|+\th^{-1}\right).
\end{align}
It follows from \eqref{ad-De} that $\Delta<1$, which implies
\begin{align}\label{Fan-add3}
(1+\th)^{ \ff T {\Delta} +2}&\leq\exp\left\{\ff {T+2} {\Delta} \log(1+\th)\right\}\cr
&\leq  \exp\left\{\ff {2(T+2)Ce^{CT}} {\th} \log(1+\th)\right\}\cr
&\leq \exp\left\{ 2(T+2)Ce^{CT}\right\}.
\end{align}
By \eqref{ad-De} again, we get
\begin{align}\label{Fan-add4}
\th^{-1}=\ff {\left(2Ce^{CT}{\|B^H\|_\be}\right)^{\ff 1 {\be}}} {2Ce^{CT}}\leq \left(2Ce^{CT}\right)^{\ff {1-\be} {\be}} {\|B^H\|_\be^{\ff 1 {\be}}}.
\end{align}
Hence, substituting \eqref{Fan-add3} and \eqref{Fan-add4} into \eqref{Fan-add2}, we deduce that there exists a constant $C>0$ such that
\begin{align}\label{ad-supn}
\|X\|_\infty\leq C\left(1+|X(0)|+\|B^H\|_\be^{\ff 1 {\be}}\right),
\end{align}
which together with the Fernique theorem  yields that for any $\de<2H,~c>0$
$$\E e^{c\|X\|_\infty^{\de}}<\infty.$$

Finally, we investigate the exponential integrability of $\|X\|_\be$.

Let $0=t_0<t_1<t_2<\cdots<t_n=T$.  we first prove that
\begin{align}\label{ad-addi}
\|f\|_\be\leq n^{1-\be}\max_{0\leq k\leq n-1}\{\|f\|_{t_k , t_{k+1} ,\be }\},~f\in C^\be([0,T]),~\be\in (0,1].
\end{align}
In fact, for any $0\leq s<t\leq T$, we have
\begin{align*}
\ff {|f(t)-f(s)|} {|t-s|^{\be}}&\leq \sum_{k=0}^{n-1}\ff {|f((t_{k+1}\wedge t)\vee s)-f((t_k\wedge t)\vee s)|} {|t-s|^{\be}}\\
&\leq \sum_{k=0}^{n-1}\ff {|(t_{k+1}\wedge t)\vee s-(t_k\wedge t)\vee s|^\be} {|t-s|^{\be}}\|f\|_{(t_k\wedge t)\vee s, (t_{k+1}\wedge t)\vee s,\be }\\
&\leq \max_{0\leq k\leq n-1}\{\|f\|_{t_k , t_{k+1} ,\be }\}\sum_{k=0}^{n-1}\ff {|(t_{k+1}\wedge t)\vee s-(t_k\wedge t)\vee s|^\be} {|t-s|^{\be}}\\
&=\max_{0\leq k\leq n-1}\{\|f\|_{t_k , t_{k+1} ,\be }\}\ff n {|t-s|^{\be}}\sum_{k=0}^{n-1}\ff {|(t_{k+1}\wedge t)\vee s-(t_k\wedge t)\vee s|^\be} {n}\\
&\leq\max_{0\leq k\leq n-1}\{\|f\|_{t_k , t_{k+1} ,\be }\}\ff n {|t-s|^{\be}}\left(\sum_{k=0}^{n-1}\ff { (t_{k+1}\wedge t)\vee s-(t_k\wedge t)\vee s } {n}\right)^{\be}\\
&= n^{1-\be}\max_{0\leq k\leq n-1}\{\|f\|_{t_k , t_{k+1} ,\be }\},
\end{align*}
where the Jenssen inequality is used in the last second inequality.\\
Applying \eqref{ad-addi} to  $X(t)$ and taking into account \eqref{ad-n-Xbe} and \eqref{ad-supn},  we have
\begin{align}\label{ad-1be}
\|X\|_\be&\leq \left({ 1+\left[\ff T {\Delta} \right]}\right)^{1-\be}\max_{0\leq k\leq \left[ \ff T {\Delta} \right]} \|X\|_{k\Delta,(k+1)\Delta,\be}\nonumber\\
&\leq \left(1+\ff T {\Delta} \right)^{1-\be}2Ce^{CT}\max_{0\leq k\leq \left[ \ff T {\Delta} \right]}\left(\|B^H\|_\be+\Delta^{1-\be}(1+|X(k\Delta)|)\right)\nonumber\\
&\leq C\left(1+\ff T {\Delta} \right)^{1-\be}\left(\|B^H\|_\be+\Delta^{1-\be}\left(1+|X(0)|+\|B^H\|_\be^{\ff 1 {\be}}\right)\right)\nonumber\\
&\leq C(1+T)^{1-\be}\left(\|B^H\|_\be \Delta^{\be-1}+1+|X(0)|+\|B^H\|_\be^{\ff 1 {\be}}\right).
\end{align}
By \eqref{ad-De}, it is clear that
$$\Delta^{\be-1}=\left(2Ce^{CT}{ \|B^H\|_\be}\right)^{{ \ff {1-\be} {\be}}},$$
which yields
$$\Delta^{\be-1}\|B^H\|_{\be}\leq \left(2Ce^{CT}\right)^{-\ff {1-\be} {\be}}{ \|B^H\|_\be^{\ff 1 {\be}}}.$$
Putting this into \eqref{ad-1be}, we obtain that
\begin{align*}
\|X\|_\be\leq C\left(1+|X(0)|+\|B^H\|_\be^{\ff 1 {\be}}\right).
\end{align*}
Therefore, for any $\de<2H$ and $c>0$, we have $\E\exp[c\|X\|_\be^{\de}]<\infty$.
\fin

\bigskip

\emph{Proof of Corollary \ref{Cor2.1}.}
By \eqref{2.2} we first obtain
\beqnn
&&K_H^{-1}\left(\int_0^\cdot h(X(r))\d r\right)(s)
=s^{H-\frac{1}{2}}D_{0+}^{H-\frac{1}{2}}\left(\cdot^{\frac{1}{2}-H} h( X(\cdot))\right)(s)\\
&=&\frac{H-\frac{1}{2}}{\Gamma(\frac{3}{2}-H)}\Bigg[\frac{1}{H-\frac{1}{2}}s^{\frac{1}{2}-H} h( X(s))
+s^{H-\frac{1}{2}} h( X(s))\int_0^s\frac{s^{\frac{1}{2}-H}-r^{\frac{1}{2}-H}}{(s-r)^{\frac{1}{2}+H}}\d r\\
&&+s^{H-\frac{1}{2}}\int_0^s\frac{ h(X(s))-h(X(r))}{(s-r)^{\frac{1}{2}+H}}r^{\frac{1}{2}-H}\d r\Bigg]
=:\frac{H-\frac{1}{2}}{\Gamma(\frac{3}{2}-H)}(I_1(s)+I_2(s)+I_3(s)).
\eeqnn
Observe that
\beqnn
\int_0^s\frac{s^{\frac{1}{2}-H}-r^{\frac{1}{2}-H}}{(s-r)^{\frac{1}{2}+H}}\d r
=\int_0^1\frac{u^{\frac{1}{2}-H}-1}{(1-u)^{\frac{1}{2}+H}}\d u\cdot s^{1-2H}=:C_0 s^{1-2H},
\eeqnn
where we make the change of variable $u=r/s, C_0$ is some constant.

(i) If $h$ satisfies \eqref{h-Holder}, then we obtain
\begin{align}\label{I1_I2}
|I_1(s)+I_2(s)|&\leq Cs^{\frac{1}{2}-H}\left(1+|X(s)|^{p+\la}\right),\nonumber\\
 \int_0^t|I_1(s)+I_2(s)|^2\d s&\leq Ct^{2(1-H)}\left(1+\|X\|_{0,t,\infty}^{2(p+\la)}\right).
\end{align}
For $I_3(s)$, let $\be<H$ such that $\la\be+\ff 1 2- H>0$. Then, we get
\begin{align*}
\int_0^s\ff {\left|  h(X(s))-h(X(r))\right|} {(s-r)^{\ff 1 2+ H}}r^{\ff 1 2-H}\d r
\leq & { K_1}\int_0^s\ff {\|X\|_{r,s,\be}^\la(s-r)^{\la\be}\left(1+|X(s)|^p+|X(r)|^p\right)} {(s-r)^{\ff 1 2+ H}r^{H-\ff 1 2}}\d r\\
\leq & { K_1}s^{\la\be+1-2H}\|X\|_{0,s,\be}^\la\left(1+\|X\|_{0,s,\infty}^p\right).
\end{align*}
Thus, it follows that
\begin{align}\label{3_I3}
|I_3(s)|\leq { K_1}s^{\la\be+\ff 1 2-H}\|X\|_{0,s,\be}^\la\left(1+\|X\|_{0,s,\infty}^p\right).
\end{align}
Let
$$\ta_n=\inf\left\{t>0~\Big|~\int_0^t\left|K_H^{-1}\left(\int_0^\cdot h(X(r))\d r\right)(s)\right|^2\d s\geq n\right\},\qquad n\in\mathbb N.$$
It is clear that under $R_{T\wedge \ta_n}\P$,
$$\tilde W_{n,t}:=W_t-\int_0^{t\wedge \ta_n} K_H^{-1}\left(\int_0^\cdot h(X(r))\d r\right)(s)\d s,~t\geq 0$$
is a Brownian motion. Then
$$\tilde B^{H}_{n,t}:=B^H_t-\int_0^{t\wedge \ta_n} h(X(r))\d r,~t\geq 0$$
is a fractional Brownian motion under $R_{T\wedge \ta_n}\P$.
Moreover, the process $X$ satisfies, for any $0\leq s\leq t\leq T$,
$$X(t)=X(s)+\int_s^t b(X(r))\d r+\int_s^t\si(X(r))\d \tilde B^H_{n,r}+\int_{s\wedge \ta_n}^{t\wedge \ta_n}\si(X(r))h(X(r))\d r.$$
Since $\si h$ has linear growth and $\tilde B^H_{n,\cdot}$ under $R_{T\wedge\ta_n}\P$ has the same distribution as $B^H_\cdot$ under $\P$,
with a slight modification of the proof of Theorem \ref{Theorem 2.1} we arrive at
$$\sup_{n}\E R_{T\wedge \ta_n} \left(\|X\|_\infty^q+\|X\|_{\be}^q\right)<\infty,\qquad q>0, \be<H.$$
Combining this with \eqref{I1_I2} and \eqref{3_I3}, we have
$$\sup_{n}\E R_{T\wedge\ta_n}\log R_{T\wedge\ta_n}
\leq C\sup_{n}\E R_{T\wedge\ta_n}\left(1+\|X\|_\infty^{2(p+\la)}+\|X\|_\be^{2\la}+\|X\|_\be^{2\la}\|X\|_\infty^{2p}
\right)<\infty.$$
Then by the Fatou lemma and the martingale convergence theorem, $\{R_t\}_{t\in [0,T]}$ is a uniformly integrable martingale and
$$\sup_{t\in [0,T]}\E R_t\log R_t<\infty.$$
As a consequence, by the Girsanov theorem we obtain that under $R_T\P$,  the process $\tilde B^H$ is a fractional Brownian motion.

(ii) By \eqref{add-n-h}, we get
\beqlb\label{Proof of Lemma 4.1-add0}
|I_1(s)+I_2(s)|\leq Cs^{\frac{1}{2}-H}\left(1+\|X\|_{0,s,\infty}^q\right).
\eeqlb
Thus, it follows from Theorem \ref{Theorem 2.1} and $q<H$ that
\beqlb\label{Proof of Lemma 4.1-1}
\E\exp\left(\int_0^T(I_1(s)+I_2(s))^2\d s\right)\leq \E e^{C\left(1+\|X\|_\infty^{2q}\right)}<\infty.
\eeqlb

Now we aim to estimate the term $I_3(s)$.

Since $H>\ff 1 2$ and $\la>\ff {2H-1} {2H}$, we can choose $\beta<H$ such that $\lambda>\frac{2H-1}{2\beta}$.
By { \eqref{add-n-h}} and \eqref{BasIne}, we derive
\begin{align}\label{h}
&\int_0^s\ff {\left|  h(X(s))-h(X(r))\right|} {(s-r)^{\ff 1 2+ H}}r^{\ff 1 2-H}\d r\nonumber\\
&\qquad\leq  { K_2}\int_0^s\ff {\left[\left(\|X\|_{r,s,\be}^\la(s-r)^{\la\be}\right)\wedge 1 \right]\left(1+|X(s)|^q+|X(r)|^q\right)} {(s-r)^{\ff 1 2+ H}r^{H-\ff 1 2}}\d r\nonumber\\
&\qquad\leq  C\left(1+\|X\|_{0,s,\infty}^q\right)\int_0^s\ff {\left(\|X\|_{0,s,\be}^\la(s-r)^{\la\be}\right)\wedge 1} {(s-r)^{\ff 1 2+ H}r^{H-\ff 1 2}}\d r\nonumber\\
&\qquad\leq   C\left(1+\|X\|_{0,s,\infty}^q\right)\int_0^s\ff {\|X\|_{0,s,\be}^\la(s-r)^{\la\be-\ff 1 2 -H}r^{\ff 1 2-H}} {1+\|X\|_{0,s,\be}^\la(s-r)^{\la\be}}\d r.
\end{align}
For the integral in the right hand of \eqref{h}, by $\lambda>\frac{2H-1}{2\beta}$ and Lemma \ref{BasIne_1} we have}
\begin{align}\label{est_int_1}
&\int_0^s\ff {\|X\|_{0,s,\be}^\la(s-r)^{\la\be-\ff 1 2 -H}r^{\ff 1 2-H}} {1+\|X\|_{0,s,\be}^\la(s-r)^{\la\be}}\d r\nonumber\\
=&\int_0^{\ff s 2}\ff {\|X\|_{0,s,\be}^\la(s-r)^{\la\be-\ff 1 2 -H}r^{\ff 1 2-H}} {1+\|X\|_{0,s,\be}^\la(s-r)^{\la\be}}\d r+\int_{\ff s 2}^s\ff {\|X\|_{0,s,\be}^\la(s-r)^{\la\be-\ff 1 2 -H}r^{\ff 1 2-H}} {1+\|X\|_{0,s,\be}^\la(s-r)^{\la\be}}\d r\nonumber\\
\leq &C\left(\int_0^{\ff s 2}\ff {\|X\|_{0,s,\be}^\la s^{\la\be-\ff 1 2 -H}r^{\ff 1 2-H}} {1+\|X\|_{0,s,\be}^\la s^{\la\be}}\d r+s^{\ff 1 2 -H}\int_{\ff s 2}^s\ff {\|X\|_{0,s,\be}^\la(s-r)^{\la\be-\ff 1 2 -H}} {1+\|X\|_{0,s,\be}^\la(s-r)^{\la\be}}\d r\right)\nonumber\\
\leq &{ C\left(s^{-\ff 1 2 -H}\int_0^{\ff s 2}r^{\ff 1 2-H}\d r+s^{\ff 1 2 -H}\int_0^{\ff s 2}\ff {\|X\|_{0,s,\be}^\la u^{\la\be-\ff 1 2 -H}} {1+\|X\|_{0,s,\be}^\la u^{\la\be}}\d u\right)}\nonumber\\
\leq & C\left(s^{1-2H}+s^{\ff 1 2-H}\|X\|_{0,s,\be}^{\ff {H-\ff 1 2} {\be}}\right).
\end{align}
Consequently, we have
\beqlb\label{4.1-0}
|I_3(s)|&\leq&C\left(s^{\ff 1 2-H}+\|X\|_{0,s,\be}^{\ff {H-\ff 1 2} {\be}}\right)\left(1+\|X\|_{0,s,\infty}^q\right)\nonumber\\
&= & C\left(s^{\ff 1 2-H}+s^{\ff 1 2-H}\|X\|_{0,s,\infty}^q+\|X\|_{0,s,\be}^{\ff {H-\ff 1 2} {\be}}+\|X\|_{0,s,\be}^{\ff {H-\ff 1 2} {\be}}\|X\|_{0,s,\infty}^q\right),
\eeqlb
which yields
\begin{align*}
\int_0^T|I_3(s)|^2\d s\leq C\left(1+\|X\|_\infty^{2q}+\|X\|_{{ \be}}^{\ff {2H-1} {\be}}+\|X\|_{ \be}^{\ff {2H-1} {\be}}\|X\|_{ \infty}^{2q}\right).
\end{align*}
Since
$$\|X\|_\infty\leq T^\be\|X\|_\be+|X(0)|,$$
we derive by { the Young inequality} that
\begin{align}\label{add-n-I3}
\int_0^T|I_3(s)|^2\d s\leq C\left(1+\|X\|_\be^{\ff {2H-1} {\be}+2q}\right).
\end{align}
Note that, by $q<H+\ff 1 {2H}-1$ we get
$$\ff {2H-1} {2(H-q)}<\ff {2H-1} {2-\ff 1 {H}}=H,$$
which allows to choose $\be$ such that
$$H>\be>\ff {2H-1} {2(H-q)}\vee\ff {2H-1} {2\la}.$$
Then we have
$$\ff {2H-1} {\be}+2q< 2 H.$$
{Consequently, combining this with \eqref{add-n-I3} and applying Theorem  \ref{Theorem 2.1} we arrive at}
\begin{align*}
\E\exp\left(\int_0^T|I_3(s)|^2\d s\right)\leq \E\exp\left[{C\left(1+\|X\|^{\ff {2H-1} {\be}+2q}_\be\right)}\right]<\infty.
\end{align*}
{Therefore, the assertion follows from this and \eqref{Proof of Lemma 4.1-1}.}
\fin

\section{Absolute continuity for fractional SDEs}

In this part, we aim to study the existence of density of the solution to  \eqref{equ_3_1} with irregular drift.
It is noticeable that we do not also impose nondegeneracy conditions on the coefficient $\sigma$.
In this context, we will invoke a method for densities, first introduced in \cite{Fournier&Printems10a} and then extended in \cite{Ro18},
which yields an easy way to prove the existence of a density.
Here, the estimates obtained in Theorem \ref{Theorem 2.1} will play an important role.

Let us then first recall the definition of the Besov space used in this section,
which is given in terms of  difference.
For any $h\in\R^d$, let $\Delta_h$ be the difference operator w.r.t. $h$ and $\Delta^m_h$ the $m$th order difference operator:
$$\Delta_h f(x)=f(x+h)-f(x),\qquad \Delta^m_h f(x)=\Delta_h\left(\Delta_h^{m-1}f\right)(x).$$
For $m\in\mathbb{N}$ and $0<\al<m$, Let $\mathscr{C}_b^\al(\R^d)$  be the Zygmund space of order $\al$ defined as the closure of bounded smooth functions w.r.t. the norm:
$$\|f\|_{\mathscr{C}_b^\al} =\|f\|_\infty+\sup_{|h|\leq 1}\ff {\|\Delta_h^mf\|_\infty} {|h|^\al},$$
and $\mathcal{B}_{1,\infty}^\al$ the Besov space of order $(1,\infty,\al)$:
$$\mathcal{B}_{1,\infty}^\al(\R^d)=\left\{f\in L^1(\R^d)~\Big|~\|f\|_{L^1}+\sup_{|h|\leq 1}\ff {\|\Delta_h^mf\|_{L^1}} {|h|^\al}<\infty\right\}.$$
The definitions above are independent of $m$,
and we can, in greater generality, define Besov space $\mathcal{B}_{p,q}^\al(\R^d)$ with $1\leq p,q\leq\infty$.
In particular, there holds $\mathscr{C}_b^\al(\R^d)=\mathcal{B}_{\infty,\infty}^\al(\R^d)$.
For more details, one can refer to \cite[Appendix A]{Ro18} and \cite{Trie1,Trie2}.

Since we do not assume that $\si$ is non-degenerate,
we will show that the distribution of $X(T)$ admits a density on $D_{\si}: =\{z\in\R^d:\sigma(z)\mbox{~is invertible}\}$
which is an open subset of $\R^d$ due to the continuity of $\si$,
and then the density belongs to some function spaces.
To this end, we
denote by $W^{s,1}(\R^d)$ (resp. $W^{s,1}_{loc}(D_{\si})$) the (resp. local) Sobolev space with $s$-order fractional derivative on $\R^d$ (resp. $D_\si$), and set
\begin{align*}
\mathcal{B}_{1,\infty,loc}^\al (D_\si)=\Big\{f~\Big|~&\mbox{for any open set}~\mathcal{O}~\mbox{with its closure}~\bar{\mathcal{O}}\subset D_\si,\\~&\mbox{there exists a}~g\in \mathcal{B}_{1,\infty}^\al(\R^d)~\mbox{such that}~g|_{\mathcal{O}}=f \Big\}.
\end{align*}

Our main result reads as follow.
\btheorem\label{Theorem 3.1}
Assume (A).

(1) Then for any $T>0$,
the law of $X(T)$ admits a density on the set $D_\si$. Moreover, the density belongs to $\mathcal{B}_{1,\infty,loc}^s(D_\si)$ for all $s<\ff {1-H} H$,
and then is in $W^{s,1}_{loc}(D_{\si})$ and $L^p_{loc}(D_\si)$ for any $s<\ff {1-H} H$ and $1\leq p<\ff {Hd} {H(d+1)-1}$.

(2) If there exists $\la>0$ such that $\si(x)\si^*(x)\geq \la$ holds for each  $x\in\R^d$,
then the law of $X(T)$ admits a density on $\R^d$ which is in $\mathcal{B}_{1,\infty}^s(\R^d)$ for any $s<\ff {1-H} H$.
Furthermore, the density is in $W^{s,1}(\R^d)$ and $L^p(\R^d)$ for any $s<\ff {1-H} H$ and $1\leq p<\ff {Hd} {H(d+1)-1}$.
\etheorem

In order to verify this theorem, some preliminary estimates are necessary.
For $\epsilon\in (0,T)$,
set
\beqnn
Y(\epsilon)&=&X(T-\epsilon)+\int_{T-\epsilon}^T\sigma(X(T-\epsilon))\d B_t^H\\
&=&X(T-\epsilon)+\sigma(X(T-\epsilon))(B_T^H-B_{T-\epsilon}^H).
\eeqnn

\blemma\label{Lemma 3.l}
For all $u\in\R^d$, there holds
\beqnn
\E\left(\e^{i\<u,Y(\epsilon)\>}|\mathcal {F}_{T-\epsilon}\right)=\exp\left\{i\<u,x\>-\ff 1 2 |\si^*(y)u|^2\int_{T-\ep}^TK_H^2(T,s)\d s\right\}\Big|_{x=\xi,~y=\et},
\eeqnn
where
$$\xi=X(T-\ep)+\si(X(T-\ep))\int_0^{T-\ep}\left(K_H(T,s)-K_H(T-\ep,s)\right)\d W_s,\ \ \et=X(T-\ep).$$
This implies that under $\P(\cdot|\mathcal {F}_{T-\epsilon}), Y(\epsilon)$ is a Gaussian random variable with the mean $\xi$ and the covariance matrix
\begin{align}\label{cov}
\mathbf{Cov}_\ep(\et)\equiv \left(\int_{T-\ep}^TK_H^2(T,s)\d s\right)\si(\et)\si^*(\et).
\end{align}
That is, for any $B\in\mathscr{B}(\R^d)$, there exists
\begin{align}
\P\left(Y(\ep)\in B\right)=\E\left(\mu_{\xi,\mathbf{Cov}_\ep(\et)}(B)\right),
\end{align}
where for $\P$-a.s. $\om\in\Omega$, the measure $\mu_{\xi(\om),\mathbf{Cov}_\ep(\et(\om))}$ is  a Gaussian measure with the  mean $\xi(\om)$ and the covariance matrix $\mathbf{Cov}_\ep(\et(\om))$.
\elemma

\emph{Proof.}
By the integral representation of $B_\cdot^H$, we get
\beqnn
B_T^H-B_{T-\epsilon}^H=\int_0^{T-\epsilon}\left(K_H(T,s)-K_H(T-\epsilon,s)\right)\d W_s+\int_{T-\epsilon}^TK_H(T,s)\d W_s.
\eeqnn
Consequently, we have
\beqlb\label{Pf-Lemma 3.l-1}
Y(\epsilon)=\xi+\si(X(T-\ep))\int_{T-\epsilon}^TK_H(T,s)\d W_s,
\eeqlb
where
$$\xi=X(T-\ep)+\si(X(T-\ep))\int_0^{T-\ep}\left(K_H(T,s)-K_H(T-\ep,s)\right)\d W_s$$
is $\mathcal {F}_{T-\epsilon}$-measurable.\\
Observing that $\int_{T-\epsilon}^TK_H(T,s)\d W_s$ is a Gaussian random variable independent of $\mathcal {F}_{T-\epsilon}$
whose covariance matrix  is $\left(\int_{T-\epsilon}^TK_H^2(T,s)\d s\right) I$, in which $I$ is the identity matrix on $\R^d$,
we obtain, for each $u\in\R^d$,
\beqnn
&&\E\left(\exp\left\{i\<u,\si(X(T-\ep))\int_{T-\ep}^T K_H(T,s)\d W_s\>\right\}\Big|\mathcal {F}_{T-\epsilon}\right)\\
&=&\E\left(\exp\left\{i\<u,\si(y)\int_{T-\ep}^T K_H(T,s)\d W_s\>\right\}\right)\Big|_{y=\et}\\
&=&\exp\left\{-\frac{|\si^*(\et) u|^2\int_{T-\epsilon}^TK_H^2(T,s)\d s}{2}\right\}.
\eeqnn
Then, combining this with \eqref{Pf-Lemma 3.l-1} yields
\beqlb\label{Pf-Lemma 3.l-2}
&&\E\left(\e^{i\<u,Y(\epsilon)\>}|\mathcal {F}_{T-\epsilon}\right)\nonumber\\
&=& \exp\left\{i\<u,\xi\>\right\}\E\left(\exp\left\{i\<u,\si(X(T-\ep))\int_{T-\ep}^T K_H(T,s)\d W_s\>\right\}|\mathcal {F}_{T-\epsilon}\right)\nonumber\\
&=&\exp\left\{i\<u,\xi\>-\frac{|\si^*(\et) u|^2\int_{T-\epsilon}^TK_H^2(T,s)\d s}{2}\right\}.
\eeqlb
Finally, it is easy to see that by \eqref{Pf-Lemma 3.l-2} the other assertions hold true.
\fin

\blemma\label{Lemma 3.2}
Assume (A). Then there holds with $\be\in[\ff 1 {1+\ga},H)$
\beqnn
\E|X(T)-Y(\epsilon)|\leq C\left(\epsilon^{(\gamma+1)\beta}+\epsilon\right).
\eeqnn
\elemma

\emph{Proof.}
Similar to \eqref{Lemma 3.2-4} and \eqref{Lemma 3.2-5}, we have
\beqnn
|X(T)-Y(\epsilon)|&=&\left|\int_{T-\epsilon}^T\left(\sigma(X(r))-\sigma(X(T-\epsilon))\right)\d B^H_r+\int_{T-\epsilon}^Tb(X(r))\d r\right|\\
&\leq&C\|B^H\|_\beta\Bigg(\int_{T-\epsilon}^T|\sigma(X(r))-\sigma(X(T-\epsilon))|(T-r)^{\alpha+\beta-1}(r-(T-\epsilon))^{-\alpha}\d r\\
&&+\|\sigma\|_\gamma\|X\|_{T-\epsilon,T,\beta}^\gamma\epsilon^{(\gamma+1)\beta}\Bigg)+C(1+\|X\|_\infty)\epsilon\\
&\leq&C\|B^H\|_\beta\|\sigma\|_\gamma\|X\|_\beta^\gamma\epsilon^{(\gamma+1)\beta}+C(1+\|X\|_\infty)\epsilon.
\eeqnn
which, together with Theorem \ref{Theorem 2.1}, leads to the desired assertion.
\fin

Besides, we shall recall a result due to \cite{Ro18},
in which the sufficient conditions for  a finite measure  admitting a density are established.
\blemma\label{Lemma 3.3}
Let $\mu$ be a finite measure on $\R^d$. If there exist  $m\in\mathbb{N}$,  $s>0,\al>0$ with $\al<s<m$, and a constant $K>0$ such that for each $\ph\in \mathscr{C}_b^\al(\R^d)$ and $h\in\R^d$ with $|h|\leq 1$,
\begin{align}\label{smooth-lem}
\left|\int_{\R^d} \Delta_h^m\ph(x) \mu(\d x)\right|\leq K|h|^s \|\ph\|_{\mathscr{C}_b^\al},
\end{align}
then  $\mu$ has a density w.r.t. the Lebesgue measure on $\R^d$.  Moreover, $\ff {\d \mu} {\d x}\in \mathcal{B}_{1,\infty}^{s-\al}(\R^d)$.
\elemma

The following criterion allows us to prove the absolute continuity of a probability measure by the localization argument.
The proof is identical to the one proposed in \cite[Lemma 1.2]{Fournier&Printems10a} and so we omit it here.
\blemma\label{Lemma 3.4}
For a continuous function $\varrho:\R^d\rightarrow\R_+$, let $D_\de=\left\{x~|~\varrho(x)\leq \de\right\}$ and $$h_\de(x)=\left(\inf\{|x-z|~|~z\in D_\de\}\right)\wedge \de,$$
where if $D_\de=\emptyset$, we assume $\inf\{|x-z|~|~z\in D_\de\}=\infty$.
Then $h_\delta:\R^d\rightarrow[0,\de]$ is a function vanishing on $D_\de$, positive on $\R^d-D_\de$ and globally Lipschitz continuous with Lipschitz constant 1. For a probability measure $\mu$ on $\R^d$,
if for each $\delta>0$, the measure $\mu_\delta(\d z)=h_\delta(z)\mu(\d z)$ has a density, then $\mu$ has a density on $\{z\in\R^d:\varrho(z)>0\}$.
\elemma
\bremark\label{Remark 3.x}
Compared with \cite[Lemma 1.2]{Fournier&Printems10a},
the Lipschitz continuity of the density $h_\de=\ff {\d \mu_\de} {\d \mu}$ in Lemma \ref{Lemma 3.4} is independent of the function $\varrho$.
\eremark

Now, we can go back to the proof of Theorem \ref{Theorem 3.1}.

\emph{Proof of Theorem \ref{Theorem 3.1}.}
Without loss of generality, we only prove the result holds true for $X(T)$.
Take $\varrho(x)$ to be the minimal eigenvalue of $\sqrt{\si(x)\si^*(x)}$ in Lemma \ref{Lemma 3.4}, and
let $\mu_\delta(\d z)=h_\delta(z)\P_{X(T)}(\d z)$.
To prove $\P_{X(T)}(\d z)$ has a density on the set $\{z\in\R^d:\sigma(z)\mbox{~is invertible}\}$,
by Lemma \ref{Lemma 3.4} it suffices to show that $\mu_\delta(\d z)$ admits a density for all $\de>0$,
which will be fulfilled via Lemma \ref{Lemma 3.3}.

For $m\in\mathbb{N}$, $\al\in (0,m)-\mathbb{N}$, $h\in\R^d$ with $ |h|<1$, and $\ph\in\mathscr{C}_b^{\al}(\R^d)$,
using Theorem \ref{Theorem 2.1} and Lemma \ref{Lemma 3.l} along with Lemma \ref{Lemma 3.2} we obtain that
\begin{align}\label{new-in-0}
\left|\E h_\de(X(T))\Delta_h^m\ph(X(T))\right| & \leq  \left|\E\left(h_\de(X(T))-h_{\de}(X(T-\ep))\right)\Delta_h^m\ph(X(T))\right|\nonumber\\
&\quad  + \left|\E h_{\de}(X(T-\ep))\left(\Delta_h^m\ph(X(T))-\Delta_h^m\ph(Y(\ep))\right)\right|\nonumber\\
&\quad + \left| \E h_{\de}(X(T-\ep))\Delta_h^m\ph(Y(\ep))\right|\nonumber\\
& \leq  C_{\de} \|\ph\|_{\mathscr{C}_b^{\al}} |h|^\al \E |X(T)-X(T-\ep)|\nonumber\\
& \quad +C_{m}\|\ph\|_{\mathscr{C}_b^{\al}}|h|^{\lfloor\al\rfloor}\E |X(T)-Y(\ep)|^{\al-\lfloor\al\rfloor}\nonumber\\
&\quad  + \left|\E h_{\de}(X(T-\ep))\E\left(\Delta_h^m\ph(Y(\ep)))| \mathcal {F}_{T-\epsilon} \right)\right|\nonumber\\
& \leq C\|\ph\|_{\mathscr{C}_b^{\al}}\left(|h|^\al\ep^\be+|h|^{\lfloor\al\rfloor}\left(\ep^{(\ga+1)\be}+\ep\right)^{\al-\lfloor\al\rfloor}\right)\nonumber\\
&\quad +\left|\E h_{\de}(\et)\mu_{\xi,\mathbf{Cov}_\ep(\et)}\left(\Delta_h^m\ph)\right)\right|,
\end{align}
where $\xi,\et$ are given in Lemma \ref{Lemma 3.l}. \\
Let $\mathbf{p}_y$ be the density of the Gaussian distribution $N(0, \mathbf{Cov}_\ep(y))$.  Then for $x\in\R^d$ and $y\in\{y|~\rh(y)\geq \de\}$, we have
\begin{align}\label{gauss-delta}
\mu_{x,\mathbf{Cov_\ep}(y)}\left(\Delta_h^m\ph)\right) & =\int_{\R^d}\Delta_h^m\ph(x+z)\mathbf{p}_y(z)\d z,\nonumber\\
&  = \int_{\R^d}\ph(x+z)\Delta_{-h}^m\mathbf{p}_y(z)\d z.
\end{align}
It is clear that
\begin{align}\label{grad-gauss}
\sup_{\rh(y)\geq \de}\int_{\R^d}|\nabla^k \mathbf{p}_y(z)| \d z & \leq C_k\sup_{\rh(y)\geq \de}\|\left(\mathbf{Cov}_\ep(y)\right)^{-\ff 1 2}\|^k\nonumber\\
& \leq C_k\left(\int_{T-\ep}^TK_H^2(T,s)\d s\right)^{-\ff {k} 2}\sup_{\rh(y)\geq \de}\rh^{-k}(y)\nonumber\\
& \leq C_k\de^{-k}\left(\int_{T-\ep}^TK_H^2(T,s)\d s\right)^{-\ff {k} 2},~k\in\mathbb{N}.
\end{align}
Notice that, by \eqref{equ-KH} we have
\beqlb\label{KH-low}
\int_{T-\epsilon}^TK_H^2(T,s)\d s
&=&\frac{1}{\Gamma^2(H-1/2)}\int_{T-\epsilon}^T\left(s^{1/2-H}\int_s^Tr^{H-1/2}(r-s)^{H-3/2}\d r\right)^2\d s\nonumber\\
&\geq&\frac{1}{\Gamma^2(H-1/2)}\int_{T-\epsilon}^T\left(\int_s^T(r-s)^{H-3/2}\d r\right)^2\d s\nonumber\\
&=&C_H\epsilon^{2H},
\eeqlb
where $C_H=\frac{1}{2H[(H-1/2)\Gamma(H-1/2)]^2}$.\\
Hence, it follows from  \eqref{gauss-delta}-\eqref{KH-low}  that
\begin{align*}
\left|\mu_{x,\mathbf{Cov_\ep}(y)}\left(\Delta_h^m\ph)\right) \right|\leq C_HC_k\de^{-k}\|\ph\|_\infty \ep^{-mH}|h|^m.
\end{align*}
Substituting this into \eqref{new-in-0} and choosing $\be$ such that $(1+\ga)\be>1$, we have
\begin{align*}
\left|\E h_\de(X(T))\Delta_h^m\ph(X(T))\right|\leq C\|\ph\|_{\mathscr{C}_b^{\al}}\left(|h|^\al\ep^\be+|h|^{[\al]}\ep^{\left(\al-[\al]\right)}+\left(\ff {|h|} {\ep^{H}}\right)^{m}\right).
\end{align*}
Let $m>\ff {\al} {1-H}$ and $\ep= \ff {1} 2 T|h|^{\ff {m-[\al]} { \al-[\al]+Hm}}$.
Then $\al<\ff {(m-[\al])(\al-[\al])} {\al-[\al]+Hm}+[\al]<m$ and
\begin{align*}
\left|\E h_\de(X(T))\Delta_h^m\ph(X(T))\right|\leq C\|\ph\|_{\mathscr{C}_b^{\al}}\left(|h|^{\al+\ff {(m-[\al])\be} {\al-[\al]+Hm}}+|h|^{\ff {(m-[\al])(\al-[\al])} {\al-[\al]+Hm}+[\al]}\right).
\end{align*}
Therefore, \eqref{smooth-lem} holds for $\mu_\de$.
As a consequence, Lemma \ref{Lemma 3.3} implies that $\mu_\de$ admits a density on $\R^d$ with $\ff {\d \mu_\de} {\d x} \in \mathcal{B}_{1,\infty}^{s(m)}(\R^d)$,
where
$$s(m)=\ff {(m-[al])\be} {\al-[\al]+Hm}\wedge \ff {\left((1-H)m-[\al]\right)(\al-[\al])} {\al-[\al]+Hm}.$$
Because of
$$\lim_{m\ra\infty}s(m) =\ff {\be} {H}\wedge \ff {(\al-[\al])(1-H)} {H},~\be<H,\al>0,$$
for any $s<\ff {1-H} H$, we can choose $m,\al,\be$ such that $s(m)=s$.  Then $\ff {\d \mu_\de} {\d x} \in \mathcal{B}_{1,\infty}^{s }(\R^d)$ for all $s<\ff {1-H} H$, which implies that the density of the distribution of $\mu$ on $D_\si$ is in $\mathcal{B}_{1,\infty,loc}^s(D_\si)$ for any $s<\ff {1-H} H$.
By \cite[2.2.2/Remark 3]{Trie1} and Sobolev embedding theorem, there hold $\mathcal{B}_{1,\infty}^{s }(\R^d)\hookrightarrow W^{s,1}(\R^d)$ and $W^{s,1}(\R^d)\hookrightarrow L^p(\R^d)$ for $1\leq p\leq \ff {d} {d-s}$.
Hence, we conclude that $\ff {\d \mu_\de} {\d x}\in W^{s,1}(\R^d)$ and $\ff {\d \mu_\de} {\d x}\in L^p(\R^d)$ for all $0<s<\ff {1-H} H$ and $1\leq p<\ff {Hd} {H(d+1)-1}$. Consequently,   the density of the distribution of $\mu$ on $D_\si$ belongs to $W^{s,1}_{loc}(D_\si)$ and $L^p_{loc}(D_\si)$ for $0<s<\ff {1-H} H$ and $1\leq p<\ff {Hd} {H(d+1)-1}$.

If there exists $\la>0$ such that $\si(x)\si^*(x)\geq \la$ for all $x\in\R^d$, then $\rh(x)\geq \sqrt{\la}$. Let $\de<\sqrt{\la}$. Then $\rh_\de \equiv \de$. It is clear that our claim holds true for this case.
\fin

\bigskip

\textbf{Acknowledgement}
The research of X. Fan was supported in part by the National Natural Science Foundation of China (Grant No. 11501009, 11871076). The research of the second author was supported in part by the National Natural Science Foundation of China (Grant No. 11771326, 11901604).

\end{document}